\documentclass[reqno]{amsart}
\usepackage{graphicx,color}
\usepackage{amsmath,amsfonts,amsbsy,url,enumerate}
\newcommand\Z{\mathbf Z}
\newcommand\R{\mathbf R}
\newcommand\Q{\mathbf Q}
\newcommand\round[2][m]{\lfloor #2\rceil_{#1}}
\newcommand{\st}{\mathrm{s.t.}}
\let\ve=\mathbf  
\newcommand\vealpha{{\boldsymbol{\alpha}}}
\newcommand\vebeta{{\boldsymbol{\beta}}}

\usepackage{ifpdf}

\usepackage{amsthm}
\newtheorem{theorem}{Theorem}%
\makeatletter
\newtheorem{lemma}{Lemma}
\renewcommand*{\c@lemma}{\c@theorem}
\renewcommand*{\p@lemma}{\p@theorem}

\renewcommand*{\c@conjecture}{\c@theorem}
\renewcommand*{\p@conjecture}{\p@theorem}

\renewcommand*{\c@proposition}{\c@theorem}
\renewcommand*{\p@proposition}{\p@theorem}

\newtheorem{corollary}{Corollary}
\renewcommand*{\c@corollary}{\c@theorem}
\renewcommand*{\p@corollary}{\p@theorem}

\renewcommand*{\c@observation}{\c@theorem}
\renewcommand*{\p@observation}{\p@theorem}

\theoremstyle{definition}

\renewcommand*{\c@problem}{\c@theorem}
\renewcommand*{\p@problem}{\p@theorem}

\renewcommand*{\c@definition}{\c@theorem}
\renewcommand*{\p@definition}{\p@theorem}

\renewcommand*{\c@remark}{\c@theorem}
\renewcommand*{\p@remark}{\p@theorem}

\newtheorem{example}{Example}
\renewcommand*{\c@example}{\c@theorem}
\renewcommand*{\p@example}{\p@theorem}

\newtheorem{algorithm}{Algorithm}
\renewcommand*{\c@algorithm}{\c@theorem}
\renewcommand*{\p@algorithm}{\p@theorem}

\makeatother

\usepackage{hyperref}

\begin{document}



\title[FPTAS for mixed-integer polynomial optimization in fixed dimension]{FPTAS 
  for optimizing polynomials over the mixed-integer points of polytopes\\
  in fixed dimension}

\author[J.A.~De Loera]{Jes\'us~A.~De Loera}
\address{J.A.~De Loera:
  University of California, Dept.~of Mathematics,
  Davis CA 95616, USA}
\email{deloera@math.ucdavis.edu}
\thanks{%
    A conference version of this article, containing a part of the results
    presented here, appeared in
    \emph{Proceedings of the 17th Annual {ACM-SIAM} Symposium on Discrete
      Algorithms, {M}iami, {FL}, January 22--24, 2006, pp. 743--748}.}
\thanks{The first author gratefully acknowledges support from NSF grant DMS-0608785, a
    2003 UC-Davis Chancellor's fellow award, the Alexander von Humboldt
    foundation, and IMO Magdeburg.}

\author[R.~Hemmecke]{Raymond Hemmecke}
\address{R.~Hemmecke: Otto-von-Guericke-Universit\"at Magdeburg,
  FMA/\penalty0 IMO,
  Universit\"ats\-platz~2, 39106 Magdeburg, Germany}
\email{hemmecke@imo.math.uni-magdeburg.de}

\author[M.~K\"oppe]{Matthias K\"oppe}
\address{M.~K\"oppe: Otto-von-Guericke-Universit\"at Magdeburg, FMA/\penalty0 IMO,
  Universit\"ats\-platz~2, 39106 Magdeburg, Germany}
\email{mkoeppe@imo.math.uni-magdeburg.de}
\thanks{The third author was supported by a Feodor Lynen Research Fellowship from the
  Alexander von Humboldt Foundation.}

\author[R.~Weismantel]{Robert Weismantel}
\address{R.~Weismantel:
  Otto-von-Guericke-Universit\"at Magdeburg, FMA/\penalty0 IMO,
  Universit\"ats\-platz~2, 39106 Magdeburg, Germany}
\email{weismant@imo.math.uni-magdeburg.de}
\thanks{The second, third and fourth authors were supported
    by the European TMR network ADONET 504438.}

\date{$\relax$Revision: 1.50.2.2 $ - \ $Date: 2007/06/15 18:49:09 $ $}

\maketitle

\begin{abstract}
  We show the existence of a fully polynomial-time approximation scheme
  (FPTAS) for the problem of maximizing a non-negative polynomial over
  mixed-integer sets in convex polytopes, when the number of variables is
  fixed.  Moreover, using a weaker notion of approximation, we show the
  existence of a fully polynomial-time approximation scheme for the problem
  of maximizing or minimizing an arbitrary polynomial over
  mixed-integer sets in convex polytopes, when the number of variables is
  fixed. 
\end{abstract}

\section{Introduction}
A well-known result by H.W.\ Lenstra Jr.\ states that  
integer {\em linear} programming problems with a fixed number of variables can be
solved in polynomial time on the input size \cite{Len83}.  Likewise,
mixed integer linear programming problems with a fixed number of integer
variables can be solved in polynomial time.
It is a natural question to ask what is the computational complexity, when the number
of variables (or the number of integer variables) is fixed, of the {\em
  non-linear} mixed integer problem%
\begin{subequations}
  \label{eq:mixed-problem}
  \begin{align}
    \max\quad & f\rlap{$(x_1,\dots,x_{d_1}, z_1,\dots,z_{d_2})$} \\
    \st\quad & A\ve x + B\ve z \leq \ve b \label{eq:mixed-problem-a}\\
    & x_i \in \R && \text{for $i = 1, \dots, d_1$}, \\
    & z_i \in \Z && \text{for $i = 1, \dots, d_2$}, \label{eq:mixed-problem-b}
  \end{align}
\end{subequations}
where $f$ is a polynomial function of maximum total degree~$D$ with
rational coefficients, and $A\in\Z^{p\times d_1}$, $B\in\Z^{p\times
  d_2}$, $\ve b\in\Z^p$.  We are interested in general polynomial objective
functions~$f$ \emph{without any convexity assumptions}.

Throughout the paper we assume that the inequality system  $A\ve x + B\ve z \leq \ve
b$ describes a convex polytope, i.e., a \emph{bounded} polyhedron, which we
denote by~$P$.  The reason for this restriction are fundamental
noncomputability results for problems involving polynomials and integer
variables.  Indeed, when we permit 
unbounded feasible regions, there cannot exist any
algorithm to decide whether there exists a feasible solution
to~\eqref{eq:mixed-problem} with $f(\ve x, 
\ve z) \geq \alpha$ (for a prescribed bound~$\alpha$), ruling out the
existence of an optimization algorithm or any approximation scheme. 
This is due to the negative answer to Hilbert's tenth
problem by Matiyasevich \cite{matiyasevich-1970,matiyasevich-1993}.  Due to
Jones' strengthening of this negative result \cite{jones-1982}, there also
cannot exist any such algorithm for the cases of unbounded feasible regions
for any fixed number of integer variables $d_2 \geq 10$. 

For the purpose of complexity analysis, we assume
that the data $A$, $B$, and $\ve b$ are given by the binary encoding
scheme, and that the objective function~$f$ is given as a list of monomials, where
the coefficients are encoded using the binary encoding scheme and the exponent
vectors are encoded using the unary encoding scheme.  In other words, the
running times are permitted to grow polynomially not only in the binary
encoding of all the problem data, but also in the maximum total
degree~$D$ of the objective function~$f$.

It is well-known that pure continuous polynomial optimization over polytopes
($d_2=0$) in varying dimension is NP-hard and that a fully polynomial time
approximation scheme ({\small FPTAS}) is not possible (unless
$\mathrm{P}=\mathrm{NP}$). Indeed the max-cut problem can be modeled as
minimizing a quadratic form over the cube $[-1,1]^d$.
H\aa{}stad~\cite{Hastad:inapprox97} proved that the max-cut problem cannot be
approximated to a ratio better than 1.0625 (unless $\mathrm{P}=\mathrm{NP}$).
This excludes the possibility of a polynomial time approximation scheme
for~\eqref{eq:mixed-problem} in varying dimension, even when the number of
integer variables is fixed.

On the other hand, pure continuous polynomial optimization problems over
polytopes ($d_2=0$) can be solved in polynomial time when the dimension $d_1$
is fixed.  This follows from a much more general result on the computational
complexity of approximating the solutions to general algebraic formulae over
the real numbers by Renegar~\cite{Renegar:1992:Approximating}; see also
\cite{Renegar:1992:CCGa,Renegar:1992:CCGb,Renegar:1992:CCGc}.

However, when we permit integer variables ($d_2 > 0$), 
it turns out that, even for fixed dimension $d_1 + d_2 = 2$ and objective
functions~$f$ of maximum total degree $D=4$, problem~\eqref{eq:mixed-problem}
is an NP-hard problem~%
\cite{deloera-hemmecke-koeppe-weismantel:intpoly-fixeddim}. 
Thus the best we can hope for, even when the number of both the continuous and
the integer variables is fixed, is an approximation
result. This paper presents the best possible such result:

\begin{theorem}[Fully polynomial-time approximation scheme]
  \label{th:mipo-fptas}
  Let the dimension $d=d_1+d_2$ be fixed.  
  \begin{enumerate}[\rm(a)]
  \item There exists a fully
  polynomial time approximation scheme ({\small FPTAS}) for the
  optimization problem~\eqref{eq:mixed-problem} for all polynomial
  functions $f(x_1,\dots,x_{d_1},\allowbreak z_1,\dots,z_{d_2})$ with rational
  coefficients  
  that are non-negative on the feasible
  region~\textup{(\ref{eq:mixed-problem-a}--\ref{eq:mixed-problem-b})}.
    
\item Moreover, the restriction to non-negative polynomials is necessary, as there
  does not even exist a polynomial time approximation scheme ({\small PTAS})
  for the maximization of {\em arbitrary} 
  polynomials over mixed-integer sets in polytopes, even for fixed dimension
  $d\geq2$, unless $\mathrm{P}=\mathrm{NP}$.
\end{enumerate}
\end{theorem}

The proof of \autoref{th:mipo-fptas} is presented in 
\autoref{section:mainthm}.  As we will see, \autoref{th:mipo-fptas} is
a non-trivial consequence of the existence of {\small FPTAS} for the
problem of maximizing a non-negative polynomial with integer
coefficients over the lattice points of a convex rational polytope.
That such {\small FPTAS} indeed exist was recently settled in
our paper~\cite{deloera-hemmecke-koeppe-weismantel:intpoly-fixeddim}.  
The knowledge of
paper~\cite{deloera-hemmecke-koeppe-weismantel:intpoly-fixeddim} is not
necessary to understand this paper but,
for convenience of the reader, we include a short summary in \autoref{section:integer-ftpas}.
Our
arguments, however, are independent of which {\small FPTAS} is used in the integral case. 

Our main approach is to use grid refinement in order to approximate the
mixed-integer optimal value via auxiliary pure integer problems.
One of the difficulties on constructing
approximations is the fact that not every sequence of grids whose
widths converge to zero leads to a convergent sequence of optimal
solutions of grid optimization problems.  This difficulty is addressed in 
\autoref{section:qualityofgrid}. 
In \autoref{section:constant} we develop techniques for bounding
differences of polynomial function values.
Section~\ref{section:mainthm} contains the proof of
\autoref{th:mipo-fptas}. 

Finally, in \autoref{section:weak}, we study a different notion of approximation.
The usual definition of an {\small FPTAS} uses the notion of
$\epsilon$-approximation that is common when considering combinatorial
optimization problems, where the 
approximation error is compared to the optimal solution value, 
\begin{equation}
  \label{eq:epsilon-approx}
  \bigl| f(\ve x_\epsilon, \ve z_\epsilon) - f(\ve x_{\max}, \ve z_{\max}) \bigr| 
  \leq \epsilon f(\ve x_{\max}, \ve z_{\max}), 
\end{equation}
where $(\ve x_\epsilon, \ve z_\epsilon)$ denotes an approximate solution and
$(\ve x_{\max}, \ve z_{\max})$ denotes a maximizer of the objective function.
In  \autoref{section:weak}, we now compare the approximation error to the
\emph{range} of the objective function 
on the feasible region,
\begin{equation}
  \label{eq:deklerk-approx-1}
  \bigl| f(\ve x_\epsilon, \ve z_\epsilon) - f(\ve x_{\max}, \ve z_{\max}) \bigr| 
  \leq \epsilon \bigl| f(\ve x_{\max}, \ve z_{\max}) - f(\ve x_{\min},\ve z_{\min}) \bigr|,
\end{equation}
where additionally $(\ve x_{\min},\ve z_{\min})$ denotes a \emph{minimizer} of
the objective function on the feasible region.
This notion of approximation was proposed by various authors 
\cite{vavasis-1993,bellare-rogaway-1993,deklerk-laurent-parillo:ptas-polynomial-simplex}. It enables us to study objective functions that are
not restricted to be non-negative on the feasible region.  
We remark that, when the objective function can take negative values on the
feasible region, \eqref{eq:deklerk-approx-1} is weaker than~\eqref{eq:epsilon-approx}.
Therefore \autoref{th:mipo-fptas}\,(b) does not rule out the existence of an
{\small FPTAS} with respect to this notion of approximation.  
Indeed we prove:
\begin{theorem}[Fully polynomial-time weak-approximation scheme]
  \label{th:mipo-wfptas}
  Let the dimension $d=d_1+d_2$ be fixed.
  Let $f$ be an arbitrary polynomial
  function with rational coefficients and maximum total degree~$D$, and 
  let $P\subset \R^d$ be a rational convex polytope.
  \begin{enumerate}[\rm(a)]
  \item In time polynomial in the input size and~$D$, it is
    possible to decide whether $f$ is constant on~$P\cap \bigl(\R^{d_1}\times \Z^{d_2}\bigr)$.  
  \item In time polynomial in
    the input size, $D$, and $\frac1\epsilon$ it is
    possible to compute a solution $(\ve x_\epsilon, \ve z_\epsilon) \in
    P\cap\bigl(\R^{d_1}\times \Z^{d_2}\bigr)$ with 
    \begin{displaymath}
      \bigl|f(\ve
      x_\epsilon, \ve z_\epsilon) - f(\ve x_{\max}, \ve z_{\max})\bigr| 
      \leq \epsilon \bigl| f(\ve x_{\max},\ve z_{\max}) - f(\ve x_{\min}, \ve z_{\min}) \bigr|.
    \end{displaymath}
  \end{enumerate}
\end{theorem}

\paragraph{Notation.}
As usual, we denote by $\Q[x_1,\dots,x_{d_1},\allowbreak z_1,\dots,z_{d_2}]$
the ring of multivariate polynomials with rational coefficients.  
For writing multivariate polynomials, we frequently use the
\emph{multi-exponent notation}, $\ve z^\vealpha = z_1^{\alpha_1} \cdots
z_d^{\alpha_d}$. 

\section{An FPTAS for the integer case}
\label{section:integer-ftpas}

The first fully polynomial-time approximation scheme for the integer case
appeared in our paper~\cite{deloera-hemmecke-koeppe-weismantel:intpoly-fixeddim}. 
It is based on Alexander Barvinok's theory for encoding all the lattice points
of a polyhedron in terms of short rational functions 
\cite{bar,BarviPom}.  The set $P\cap\Z^d$ is represented
by a Laurent polynomial $g_P(\ve z)=\sum_{\vealpha \in P\cap\Z^d}
\ve z^\vealpha$.  
From Barvinok's theory this exponentially-large sum of
monomials $g_P(\ve z)$ can be written as a polynomial-size sum of
rational functions (assuming the dimension $d$ is fixed) of the form:
\begin{equation} \label{eq:aa}
g_P(\ve z) = \sum_{i\in I} {E_i \frac{\ve z^{\ve u_i}} {\prod_{j=1}^d
(1-\ve z^{\ve v_{ij}})}},
\end{equation}
where $I$ is a polynomial-size indexing set, and where
$E_i\in\{1,-1\}$ and $\ve u_i, \ve v_{ij} \in\Z^d$ for all $i$ and $j$.
There is a polynomial-time algorithm for computing this representation
\cite{bar,BarviPom,latte2,latte1}.

By symbolically applying differential operators to the
representation~\eqref{eq:aa}, we can compute a short rational function
representation of the Laurent polynomial 
\begin{equation}
  g_{P,f}(\ve z) = \sum_{\vealpha\in P \cap \Z^d} f(\vealpha) \ve z^{\vealpha}.
\end{equation}
In fixed dimension, 
the size of the expressions occuring in the symbolic calculation can be
bounded polynomially:

\begin{lemma}[\cite{deloera-hemmecke-koeppe-weismantel:intpoly-fixeddim}, Lemma 3.1] 
  \label{operators} Let the dimension $d$ be fixed. 
  Let $g_P(\ve z)=\sum_{\vealpha\in P \cap
    \Z^d}\ve z^{\vealpha}$ be the Barvinok 
  representation of the generating function of $P\cap\Z^d$. Let
  $f\in\Z[x_1,\dots,x_d]$ be a polynomial of maximum total
  degree~$D$. We can compute, in time polynomial in $D$ and the input size, a
  Barvinok representation $g_{P,f}(\ve z)$ 
  for the generating function $\sum_{\vealpha\in P \cap \Z^d} f(\vealpha) \ve z^{\vealpha}.$
\end{lemma}

Now we present the algorithm to obtain bounds $U_k,L_k$ that reach the optimum.  
We make use of the elementary fact that, for a set
$S=\{s_1,\dots,s_r\}$ of non-negative real numbers, 
\begin{equation}
  \max\{s_1, \dots, s_r \} = \lim_{k \rightarrow \infty} \sqrt[k]{s_1^r +
    \dots + s_k^r}.
\end{equation}

\begin{algorithm}[Computation of bounds]~\smallskip
\label{Algorithm}

\noindent {\em Input:} A rational convex polytope $P \subset \R^d$, 
a polynomial objective $f \in \Z[x_1,\dots,x_d]$ of maximum total degree $D$ 
that is non-negative over $P\cap\Z^d$.\smallskip

\noindent {\em Output:} A nondecreasing sequence of lower bounds $L_k$,
and a nonincreasing sequence of upper bounds $U_k$, both reaching the maximal
function value $f^*$ of $f$ over $P\cap\Z^d$ in a finite number of steps.

\begin{enumerate}
\item  Compute a short rational function expression for
  the generating function $g_P(\ve z)=\sum_{\vealpha\in P\cap\Z^d} \ve
  z^{\vealpha}$.  Using residue techniques, compute $|P \cap
  \Z^d|=g_P(\ve 1)$ from $g_P(\ve z)$.

\item From the rational function $ g_P(\ve z)$ 
  compute the rational function representation of $g_{P,f^k}(\ve z)$ of
  $\sum_{\vealpha\in P \cap \Z^d} f^k(\vealpha) \ve z^\vealpha$ by 
  \autoref{operators}. Using residue techniques, compute 
\[
L_k:=\left\lceil{\sqrt[k]{g_{P,f^k}(\ve 1)/g_{P,f^0}(\ve1)}}\,\right\rceil\;\;\;\text{and}\;\;\;
U_k:=\left\lfloor{\sqrt[k]{g_{P,f^k}(\ve 1)}}\right\rfloor. 
\] 
\end{enumerate} 
\end{algorithm}

\begin{theorem}[\cite{deloera-hemmecke-koeppe-weismantel:intpoly-fixeddim}, Lemma 3.3 and Theorem 1.1]\label{lemma:bounds} 
  Let the dimension $d$ be fixed.  Let $P\subset\R^d$ be a rational convex polytope.
  Let $f$ be a polynomial with integer coefficients and maximum total
  degree~$D$ that is non-negative on $P\cap\Z^d$. 
  \begin{enumerate}[\rm(i)]
\item \autoref{Algorithm} computes the bounds $L_k$, $U_k$ in time polynomial in
  $k$, the input size of $P$ and $f$, and the total degree~$D$. The bounds
  satisfy the following inequality:
$$ 
U_k-L_k \leq f^* \cdot \left(\sqrt[k]{|P \cap \Z^d|}-1 \right).
$$
\item For $k=(1+1/\epsilon)\log({|P \cap \Z^d|})$ (a number bounded by a
  polynomial in the input size),
  $L_k$ is a $(1-\epsilon)$-approximation to the optimal value $f^*$ and it
  can be computed in time polynomial in the input size, the total
  degree~$D$, and $1/\epsilon$. Similarly, $U_k$ gives a
  $(1+\epsilon)$-approximation to $f^*$. 

\item With the same complexity,  by iterated bisection of $P$, we can also find
  a feasible solution $\ve x_\epsilon\in P\cap\Z^d$ with 
  \begin{displaymath}
    \bigl|f(\ve x_\epsilon) - f^*\bigr| \leq \epsilon f^*.
  \end{displaymath}
\end{enumerate}
\end{theorem}

\section{Grid approximation results}
\label{section:qualityofgrid}

An important step in the development of an {\small FPTAS} for the mixed-integer
optimization problem is the reduction of the mixed-integer
problem~\eqref{eq:mixed-problem} to an auxiliary optimization problem over a
lattice $\frac1m \Z^{d_1}\times \Z^{d_2}$.  To this end, we consider the
\emph{grid problem} with grid size~$m$,
\begin{equation}
  \label{eq:grid-problem}
  \begin{aligned}
    \max\quad & f\rlap{$(x_1,\dots,x_{d_1},z_1,\dots,z_{d_2})$} \\
    \st\quad & A\ve x + B\ve z\leq \ve b \\
    & x_i \in \tfrac 1m \Z &&\text{for $i=1,\dots,d_1$},\\
    & z_i \in \Z &&\text{for $i=1,\dots,d_2$}.\\
  \end{aligned}
\end{equation}
We can solve this problem approximately using the integer {\small FPTAS}
(\autoref{lemma:bounds}):
\begin{corollary}\label{cor:gridproblem-approximation}
  For fixed dimension $d=d_1+d_2$ there 
  exists an algorithm with running time polynomial in $\log m$, the
  encoding length of~$f$ and of~$P$, the maximum total degree~$D$ of~$f$,
  and~$\frac1\epsilon$ for computing a  
  feasible solution $(\ve x^m_\epsilon, \ve z^m_\epsilon)\in
  P\cap\bigl(\frac1m\Z^{d_1}\times \Z^{d_2}\bigr)$ to the grid
  problem~\eqref{eq:grid-problem} with an objective function~$f$ that is
  non-negative on the feasible region, with
  \begin{equation}\label{eq:gridproblem-approximation}
    f(\ve x^m_\epsilon, \ve z^m_\epsilon) \geq (1-\epsilon) f(\ve x^m,\ve z^m),
  \end{equation}
  where $(\ve x^m, \ve z^m)\in
  P\cap\bigl(\frac1m\Z^{d_1}\times \Z^{d_2}\bigr)$ is an optimal solution
  to~\eqref{eq:grid-problem}.
\end{corollary}
\begin{proof} 
  We apply \autoref{lemma:bounds} to the pure
  integer optimization problem:  
  \begin{equation} 
    \label{eq:pure-integer-problem}
    \begin{aligned}
      \max\quad & \tilde f(\ve{\tilde x}, \ve z) \\
      \st\quad  & A\ve{\tilde x} + m B\ve z \leq m\ve b \\
      & \tilde x_i \in \Z &&\text{for $i=1,\dots,d_1$}, \\
      & z_i \in \Z &&\text{for $i=1,\dots,d_2$},
    \end{aligned} 
  \end{equation}
  where $\tilde f(\ve{\tilde x},\ve z) := m^D f( \tfrac1m \ve{\tilde x}, \ve
  z)$ is a polynomial function with integer coefficients. 
  Clearly the binary encoding length of the coefficients 
  of~$\tilde f$ increases by at most $\lceil D\log m \rceil$, compared to the
  coefficients 
  of~$f$.  Likewise, the encoding length of the coefficients
  of~$mB$ and~$m\ve b$ increases by at most $\lceil \log m\rceil$.
  By Theorem 1.1
  of~\cite{deloera-hemmecke-koeppe-weismantel:intpoly-fixeddim}, there exists
  an algorithm with running time polynomial in the encoding length of~$\tilde
  f$ and of~$A\ve x + m B\ve z \leq m \ve b$, the maximum total degree~$D$,
  and~$\frac1\epsilon$ for computing a feasible solution $(\ve x^m_\epsilon, \ve z^m_\epsilon)\in
  P\cap\bigl(\frac1m\Z^{d_1}\times \Z^{d_2}\bigr)$ such that $\tilde f(\ve
  x^m_\epsilon, \ve z^m_\epsilon) \geq (1-\epsilon) \tilde f(\ve x^m,\ve z^m)$,
  which implies the estimate~\eqref{eq:gridproblem-approximation}.
\end{proof} 

One might be tempted to think that for large-enough choice of~$m$, we
immediately obtain an approximation to the mixed-integer optimum with
arbitrary precision.  However, this is not true, as the following example
demonstrates. 

\begin{example} 
  Consider the mixed-integer optimization problem 
  \begin{equation} 
    \label{eq:example-slice-not-fulldim}
    \begin{aligned}
      \max\quad& 2z - x\\
      \st\quad& z \leq 2x \\
      & z \leq 2(1-x) \\
      & x\in\R_{\geq0},\ z\in\{0,1\}, 
    \end{aligned} 
  \end{equation}
  whose feasible region consists of the point $(\frac12,1)$ and the segment
  $\{\,(x,0): x\in[0,1]\}$.  The unique optimal solution
  to~\eqref{eq:example-slice-not-fulldim} is $x=\frac12$, $z=1$.  Now consider
  the sequence of grid approximations of~\eqref{eq:example-slice-not-fulldim}
  where $x\in \frac1m \Z_{\geq0}$.  For even $m$, the unique optimal solution to the
  grid approximation is $x=\frac12$, $z=1$.  However, for odd $m$, the unique
  optimal solution is $x=0$, $z=0$.  Thus the full sequence of the optimal
  solutions to the grid approximations does not converge since it has two limit
  points; see \autoref{fig:example-slice-not-fulldim}. 
  \begin{figure}[t]
    \begin{minipage}[t]{.4\linewidth}
      \ifpdf
    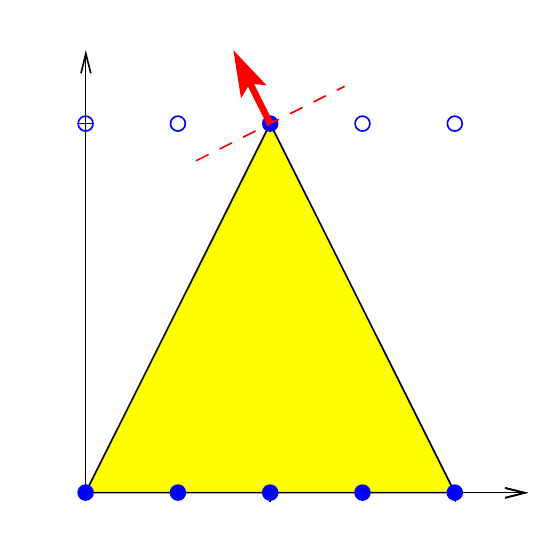
    \else
    \input{miptriangle-even.pstex_t}
    \fi
    \end{minipage}
    \qquad\qquad
    \begin{minipage}[t]{.4\linewidth}
      \ifpdf
    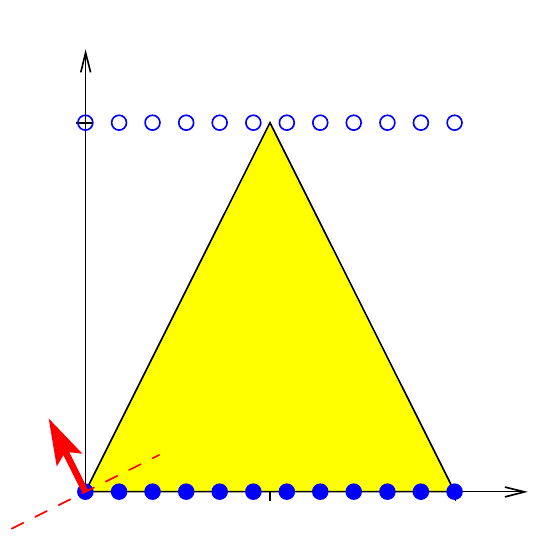
    \else
    \input{miptriangle-11.pstex_t}
    \fi
    \end{minipage}
    \caption{A sequence of optimal solutions to grid problems with two limit
      points, for even $m$ and for odd $m$}
    \label{fig:example-slice-not-fulldim}
  \end{figure}
\end{example} 

\begin{figure}[t]
  \centering
  \ifpdf
    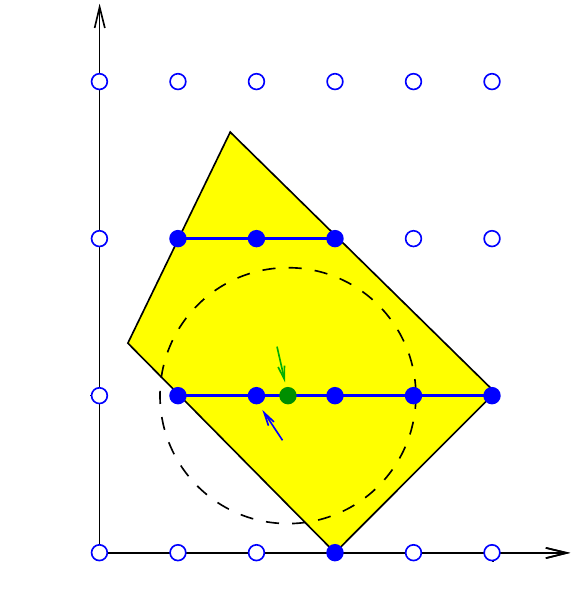
    \else
    \input{gridapprox.pstex_t}
    \fi
  \caption{The principle of grid approximation.  
    Since we can refine the grid only in the direction of the continuous
    variables, we need to construct an approximating grid point $(\ve x,\ve
    z^*)$ in the same integral slice as the target point $(\ve x^*,\ve z^*)$.}
  \label{fig:grid-approximation-lowerdim}
\end{figure}

Even though taking the limit does not work, taking the upper limit does.
More strongly, we can prove that it is possible to construct, in polynomial time, a
subsequence of finer and finer grids that contain a lattice point $(\round[\delta]{\ve x^*}, \ve z^*)$ that is arbitrarily close to the mixed-integer optimum $(\ve x^*,\ve
z^*)$.  This is the central statement of this section and a basic building
block of the approximation result.
\begin{theorem}[Grid Approximation]
  \label{lemma:grid-approximation-lowerdim}
  Let $d_1$ be fixed.
  Let $P = \{\, (\ve x,\ve z)\in\R^{d_1+d_2} : A\ve x + B\ve z \leq \ve
  b\,\}$, where $A\in\Z^{p\times d_1}$, $B\in\Z^{p\times d_2}$.  Let $M\in\R$
  be given such
  that $P \subseteq \{\, (\ve x,\ve z)\in\R^{d_1+d_2} : |x_i| \leq M \text{ for
  } i=1,\dots,d_1 \,\}$.  
  There exists a polynomial-time algorithm to compute a number $\Delta$
  such that for every $(\ve x^*,\ve z^*)\in P\cap(\R^{d_1}\times \Z^{d_2})$
  and $\delta>0$ the following property holds:
  \begin{quote} 
    Every lattice $\frac{1}{m}\Z^{d_1}$ for $m=k\Delta$ and $k\geq\frac2\delta
    d_1 M$ contains a lattice point $\round[\delta]{\ve x^*}$ such that
    $(\round[\delta]{\ve x^*}, \ve z^*) \in P\cap
    \bigl(\frac{1}{m}\Z^{d_1}\times \Z^{d_2}\bigr) $ and
    $\bigl\|\round[\delta]{\ve x^*}-\ve x^*\bigr\|_\infty \leq\delta$.
  \end{quote} 
\end{theorem}
The geometry of \autoref{lemma:grid-approximation-lowerdim} is illustrated in
\autoref{fig:grid-approximation-lowerdim}.  The notation $\round[\delta]{\ve
  x^*}$ has been chosen to suggest that the coordinates of $\ve x^*$ have been
``rounded'' to obtain a nearby lattice point.  The rounding method 
is provided by the next two lemmas;
\autoref{lemma:grid-approximation-lowerdim} follows directly from them.

\begin{lemma}[Integral Scaling Lemma]
  \label{lemma:integral-scaling}
  Let $P = \{\, (\ve x,\ve z)\in\R^{d_1+d_2} : A\ve x + B\ve z \leq \ve
  b\,\}$, where $A\in\Z^{p\times d_1}$, $B\in\Z^{p\times d_2}$.
  For fixed $d_1$, there exists a polynomial time algorithm to compute a number
  $\Delta\in\Z_{>0}$ such that for every $\ve z\in\Z^{d_2}$ the polytope
  \begin{displaymath}
    \Delta P_{\ve z} = \bigl\{\, \Delta\ve x : (\ve x,\ve z) \in P\,\bigr\}
  \end{displaymath}
  is integral, i.e., all vertices have integer coordinates.  In particular, the
  number $\Delta$ has an encoding 
  length that is bounded by a polynomial in the encoding length
  of~$P$.
\end{lemma}
\begin{proof}
  Because the dimension $d_1$ is fixed, there exist only polynomially many
  simplex bases of the inequality system $A \ve x \leq \ve b - B\ve z $, and they can be
  enumerated in polynomial time.  The determinant of each simplex basis can be
  computed in polynomial time.  Then $\Delta$ can be chosen as the least
  common multiple of all these determinants.
\end{proof}

\begin{lemma}
  \label{lemma:grid-approximation-integral}
  Let $Q\subset \R^d$ be an integral polytope.  Let $M\in\R$ be such that $Q\subseteq \{\, \ve x\in\R^d : |x_i| \leq M \text{ for
  } i=1,\dots,d \,\}$.  Let $\ve x^*\in Q$ and let $\delta>0$.  Then every lattice
  $\frac1k\Z^d$ for $k \geq \frac2\delta d M$ contains a lattice
  point $\ve x\in Q\cap \frac1k\Z^d$ with $\|\ve x-\ve x^*\|_\infty \leq\delta$.
\end{lemma} 
\begin{proof}
  By Carath\'eodory's Theorem, 
  there exist $d+1$ vertices $\ve x^0, \dots, \ve x^d\in\Z^d$ of~$Q$ and convex multipliers
  $\lambda_0,\dots,\lambda_d$ such that $\ve x^* = \sum_{i=0}^d \lambda_i \ve
  x^i$.  Let $\lambda'_i := \frac1k \lfloor k\lambda_i \rfloor \geq 0$ for
  $i=1,\dots,d$ and $\lambda'_0 := 1 - \sum_{i=1}^d \lambda_i'
  \geq0$. Moreover, we conclude $\lambda_i-\lambda'_i\leq\frac1k$ for
  $i=1,\dots,d$ and $\lambda'_0-\lambda_0=\sum_{i=1}^d
  (\lambda_i-\lambda'_i)\leq d\frac 1k$.
  Then $\ve x := \sum_{i=0}^d \lambda'_i \ve x^i \in Q\cap\frac1k\Z^d$, and
  we have
  \begin{align*}
    \|\ve x-\ve x^*\|_\infty & \leq \sum_{i=0}^d |\lambda'_i - \lambda_i| \| \ve x^i\|_\infty
    \leq 2d \frac 1k M
    \leq \delta,
  \end{align*}
  which proves the lemma.
  
\end{proof}

\section{Bounding techniques}
\label{section:constant}

Using the results of \autoref{section:qualityofgrid} we are now able to approximate the
mixed-integer optimal point by a point of a suitably fine lattice.  The
question arises how we can use the geometric distance of these two points to
estimate the difference in objective function values.  We prove \autoref{lemma:lipschitz}
that provides us with a local Lipschitz constant for the polynomial to be
maximized.

\begin{lemma}[Local Lipschitz constant]\label{lemma:lipschitz}
  Let $f$ be a polynomial in $d$ variables with maximum
  total degree~$D$. 
  Let $C$ denote the largest absolute value of a coefficient of~$f$.
  Then there exists a Lipschitz constant $L$ such that $|f(\ve x) - f(\ve y)|
  \leq L \|\ve x-\ve y\|_\infty$
  for all $|x_i|, |y_i|\leq M$.  The constant $L$ is $O(D^{d+1} CM^D)$.
\end{lemma}

\begin{proof} 
  Let 
  \begin{math}
    f(\ve x) = \sum_{\vealpha\in\mathcal D} c_{\vealpha} \ve x^\vealpha,
  \end{math}
  where $\mathcal D\subseteq\Z_{\geq0}^d$ is the set of exponent vectors of
  monomials appearing in~$f$.
  Let $r = |\mathcal D|$ be the number of monomials of~$f$.
  Then we have
  \begin{displaymath}
    |f(\ve x) - f(\ve y)| \leq \sum_{\vealpha\neq\ve0} |c_\vealpha|\, |\ve x^\vealpha - \ve y^\vealpha|.
  \end{displaymath}
  We estimate all summands separately. 
  Let $\vealpha\neq\ve0$ be an exponent vector with $n := 
  \sum_{i=1}^d \alpha_i \leq D$.  
  Let 
  \begin{displaymath}
    \vealpha = \vealpha^0 \geq \vealpha^1 \geq \dots \geq \vealpha^n = \ve 0
  \end{displaymath}
  be a decreasing chain of exponent vectors with $\vealpha^{i-1} - \vealpha^{i} = \ve
  e^{j_i}$ for $i=1,\dots,n$.  Let $\vebeta^i := \vealpha - \vealpha^i$ for $i=0,\dots,n$. 
  Then $\ve x^\vealpha-\ve y^\vealpha$ can be
  expressed as the ``telescope sum'' 
  \begin{align*}
    \ve x^\vealpha-\ve y^\vealpha
    &= \ve x^{\vealpha^0}\ve y^{\vebeta^0} - \ve x^{\vealpha^1} \ve y^{\vebeta^1} 
    + \ve x^{\vealpha^1} \ve y^{\vebeta^1} - \ve x^{\vealpha^2} \ve
    y^{\vebeta^2} +- \cdots
    - \ve x^{\vealpha^n} \ve y^{\vebeta^n} \\
    &= \sum_{i=1}^n \left( \ve x^{\vealpha^{i-1}} \ve y^{\vebeta^{i-1}} 
      - \ve x^{\vealpha^{i}} \ve y^{\vebeta^{i}} \right) \\
    &= \sum_{i=1}^n \left( (x_{j_i} - y_{j_i}) \ve x^{\vealpha^{i}} \ve y^{\vebeta^{i-1}} \right). 
  \end{align*}
  Since $ \bigl| \ve x^{\vealpha^{i}} \ve
  y^{\vebeta^{i-1}} \bigr| \leq M^{n-1}$ and $n\leq D$, we obtain 
  \begin{align*}
    | \ve x^\vealpha - \ve y^\vealpha | & \leq D \cdot \|\ve x-\ve y\|_\infty
    \cdot M^{n-1} ,
  \end{align*}
  thus
  \begin{align*}
    |f(\ve x) - f(\ve y)| \leq C r D M^{D-1} \|\ve x-\ve y\|_\infty.
  \end{align*}
  Let $L := C r D M^{D-1} $.  Now, since $r = O(D^d)$, we have $L = O(D^{d+1} CM^D)$.
\end{proof}

Moreover, in order to obtain an {\small FPTAS}, we need to put differences of
function values in relation to the maximum function value.  
To do this, we need to deal with the special case of polynomials that are constant on
the feasible region; here trivially every feasible solution is optimal.
For non-constant polynomials, we can prove a lower bound on the maximum
function value.  The technique is to bound the difference of the minimum and
the maximum function value on the mixed-integer set from below; if the
polynomial is non-constant, this implies, for a non-negative polynomial, a lower
bound on the maximum function value. We will need a simple fact about the
roots of multivariate polynomials. 

\begin{lemma}\label{lemma:toomanyroots}
  Let $f\in\Q[x_1,\dots,x_{d}]$ be a polynomial and let $D$ be the largest
  power of any variable that appears in~$f$.  Then $f=0$ if and only if $f$
  vanishes on the set $\{0,\dots, D\}^d$. 
\end{lemma}
\begin{proof}
  This is a simple consequence of the Fundamental Theorem of Algebra.  See,
  for instance, \cite[Chapter 1, \S1, Exercise 6\,b]{CoxLittleOShea92}.
\end{proof}

\begin{figure}[t]
  \centering
  \ifpdf
    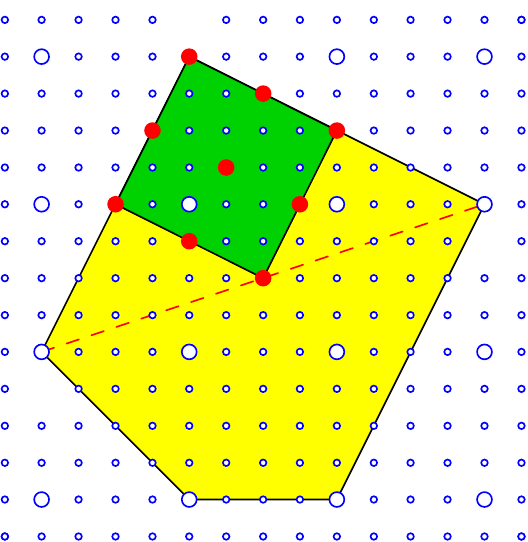
    \else
    \input{constant1-4.pstex_t}
    \fi
  \caption{The geometry of \autoref{lemma:constantness}.
    For a polynomial with maximum total degree of~$2$, we construct a
    refinement $\frac1k\Z^d$ (small circles) of the standard lattice (large circles)
    such that $P\cap\frac1k\Z^d$ contains an affine image of the set
    $\{0,1,2\}^d$ (large dots).}
  \label{fig:const}
\end{figure}

\begin{lemma}\label{lemma:constantness}
  Let $f\in\Q[x_1,\dots,x_{d}]$ be a polynomial with maximum total degree~$D$.
  Let $Q \subset \R^d$ be an integral 
  polytope of dimension $d'\leq d$.  
  Let $k \geq D\, d'$.
  Then $f$ is constant on $Q$ if and only if 
  $f$ is constant on $Q\cap \frac1k \Z^d$.  
\end{lemma}
\begin{proof}
  Let $\ve x^0\in Q\cap\Z^d$ be an arbitrary vertex of~$Q$.  There exist
  vertices $\ve x^1, \dots, \ve x^{d'}\in Q\cap\Z^d$ such that the vectors
  $\ve x^1-\ve x^0, \dots, \ve x^{d'}-\ve x^0\in\Z^d$ are linearly
  independent.  By convexity, $Q$ contains the parallelepiped 
  \begin{displaymath}
    S := 
    \Bigl\{\, \ve x^0 + \textstyle\sum_{i=1}^{d'} \lambda_i (\ve x^i-\ve x^0) :
    \text{$\lambda_i \in [0,\tfrac1{d'}]$ for $i=1,\dots, d'$} \,\Bigr\}.
  \end{displaymath}
  We consider the set 
  \begin{align*}
    S_k &= \tfrac1k \Z^d \cap S
    \supseteq 
    \Bigl\{\, \ve x^0 + \textstyle\sum_{i=1}^{d'}
    \frac{n_i}k (\ve x^i-\ve x^0) :
     \text{$n_i \in \{0,1,\dots,D\}$ for $i=1,\dots, d'$}
    \,\Bigr\};
  \end{align*}
  see \autoref{fig:const}.
  Now if there exists a $c\in\R$ with $f(\ve x) = c$ for all $\ve x\in Q\cap
  \frac1k \Z^d$, then all the points in $S_k$ are roots of the polynomial 
  $f - c$, which has only maximum total degree~$D$.  
  By \autoref{lemma:toomanyroots} (after an affine transformation), $f-c$ is
  zero on the affine hull of $S_k$; 
  hence $f$ is constant on the polytope~$Q$. 
\end{proof}

\begin{theorem}\label{lemma:mixed-constantness}
  Let $f\in\Z[x_1,\dots,x_{d_1},z_1,\dots,z_{d_2}]$.  
  Let $P$ be a rational convex polytope, and let $\Delta$
  be the number from \autoref{lemma:integral-scaling}.
  Let $m = k\Delta$ with $k \geq D\, d_1$, $k\in\Z$.
  Then $f$ is constant on the feasible region $P\cap \bigl(\R^{d_1}\times \Z^{d_2}\bigr)$ if and
  only if $f$ is constant on $P\cap \bigl( \frac1m \Z^{d_1} \times \Z^{d_2}
  \bigr)$.  If $f$ is not constant, then 
  \begin{equation}\label{eq:lower-bound-for-range}
    \bigl| f(\ve x_{\max}, \ve z_{\max})
    - f(\ve x_{\min},\ve z_{\min}) \bigr| \geq m^{-D},
  \end{equation}
  where $(\ve x_{\max}, \ve z_{\max})$ is an optimal solution to the
  maximization problem over the feasible region $P\cap \bigl(\R^{d_1}\times
  \Z^{d_2}\bigr)$ and $(\ve x_{\min}, \ve 
  z_{\min})$ is an optimal solution to the minimization problem.
\end{theorem}
\begin{proof}
  Let $f$ be constant on $P\cap \bigl( \frac1m \Z^{d_1} \times \Z^{d_2}
  \bigr)$.  For fixed integer part $\ve z\in\Z^{d_2}$, we consider the
  polytope $\Delta P_{\ve z} = \bigl\{\, \Delta\ve x :
  (\ve x,\ve z) \in P\,\bigr\}$, which is a slice of $P$ scaled to become an
  integral polytope.  By applying \autoref{lemma:constantness} with $k = (D+1)
  d$ on every polytope $\Delta P_{\ve z}$, we obtain that $f$ is constant on every slice~$P_{\ve z}$. 
  Because $f$ is also constant on the set $P\cap \bigl( \frac1m \Z^{d_1} \times \Z^{d_2}
  \bigr)$, which contains a point of every non-empty slice $P_{\ve z}$,
  it follows that $f$ is constant on $P$. 
  
  If $f$ is not constant, there exist $(\ve x^1, \ve z^1)$, $(\ve x^2, \ve
  z^2) \in P\cap\bigl( \frac1m \Z^{d_1} \times \Z^{d_2} \bigr)$ with $f(\ve x^1,\ve z^1)
  \neq f(\ve x^2,\ve z^2)$.  
  By the integrality of all coefficients of $f$, we
  obtain the estimate
  \begin{align*}
    |f(\ve x^1,\ve z^1) - f(\ve x^2,\ve z^2)| \geq m^{-D}.
  \end{align*}
  Because $(\ve x^1, \ve z^1)$, $(\ve x^2, \ve z^2)$ are both feasible solutions
  to the maximization problem and the minimization problem, this
  implies~\eqref{eq:lower-bound-for-range}. 
\end{proof}

\section{Proof of \autoref{th:mipo-fptas}} \label{section:mainthm}
 
Now we are in the position to prove the main result.

\begin{proof}[Proof of \autoref{th:mipo-fptas}]
  \emph{Part (a).}
  Let $(\ve x^*, \ve z^*)$ denote an optimal solution to the mixed-integer
  problem~\eqref{eq:mixed-problem}.  Let $\epsilon>0$.  We show that, in time
  polynomial in the input length, the maximum total degree,
  and~$\frac1\epsilon$, 
  we can compute a point 
  $(\ve x, \ve z)$ that satisfies
  \mbox{(\ref{eq:mixed-problem-a}--\ref{eq:mixed-problem-b})}
  such that 
  \begin{equation}
    \label{eq:epsilon-approximation}
    |f(\ve x, \ve z) - f(\ve x^*, \ve z^*)| \leq \epsilon f(\ve x^*, \ve z^*).
  \end{equation}
  We prove this by establishing several estimates, which are illustrated in
  \autoref{fig:estimates}. 

  \begin{figure}[t]
    \centering
    \ifpdf
    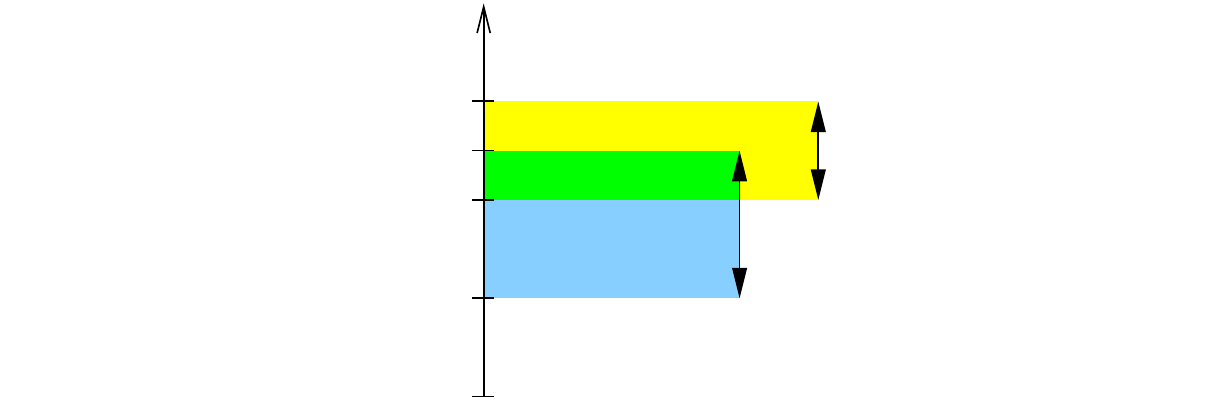
    \else
    \input{mipestimate5.pstex_t}
    \fi
    \caption{Estimates in the proof of \autoref{th:mipo-fptas}\,(a)}
    \label{fig:estimates}
  \end{figure}
  First we note that we can restrict ourselves to the case of polynomials
  with integer coefficients, simply by multiplying $f$ with the least common
  multiple of all denominators of the coefficients.
  We next establish a lower bound on $f(\ve x^*, \ve z^*)$. 
  To this end, let $\Delta$ be the integer from
  \autoref{lemma:integral-scaling}, which can be computed in
  polynomial time.  By \autoref{lemma:mixed-constantness} with $m = D\, d_1
  \Delta$, either $f$ is constant on the feasible region, or 
  \begin{equation}\label{eq:lower-bound-for-optimum}
    f(\ve x^*, \ve z^*) \geq (D\, d_1 \Delta)^{-D},
  \end{equation}
  where $D$ is the maximum total degree of~$f$.
  Now let 
  \begin{equation}
    \delta:=\frac{\epsilon}{2(D d_1 \Delta)^D L(C,D,M)}
  \end{equation}
  and let us choose the grid size
  \begin{equation}
    m := \Delta \left\lceil \frac4\epsilon (D d_1 \Delta)^D  L(C,D,M) d_1 M
    \right\rceil,
  \end{equation}
  where $L(C,D,M)$ is the Lipschitz constant from \autoref{lemma:lipschitz}.
  Then we have $m \geq \Delta \frac2\delta d_1 M$, so
  by \autoref{lemma:grid-approximation-lowerdim}, there
  is a point $(\round[\delta]{\ve x^*}, \ve z^*) \in P \cap \bigl( \frac1m
  \Z^{d_1}\times \Z^{d_2}\bigr)$  
  with $\bigl\| \round[\delta]{\ve x^*} - \ve x^* \bigr\|_\infty \leq \delta$. 
  Let $(\ve x^m, \ve z^m)$ denote an optimal solution to the grid problem~\eqref{eq:grid-problem}.  
  Because $(\round[\delta]{\ve x^*}, \ve z^*)$ is a feasible solution to the grid problem~\eqref{eq:grid-problem}, 
  we have 
  \begin{equation}
    f(\round[\delta]{\ve x^*}, \ve z^*) \leq f(\ve x^m, \ve z^m) \leq f(\ve x^*, \ve z^*).
  \end{equation}
  Now we can estimate
  \begin{align}
    \bigl| f(\ve x^*, \ve z^*) - f(\ve x^m, \ve z^m) \bigr| 
    & \leq \bigl| f(\ve x^*, \ve z^*) - f(\round[\delta]{\ve x^*}, \ve z^*) \bigr| \nonumber\\
    & \leq L(C,D,M) \, \bigl\| \ve x^* - \round[\delta]{\ve x^*} \bigr\|_\infty \nonumber\\
    & \leq L(C,D,M) \, \delta \nonumber\\
    & = \frac\epsilon2 {(D\, d_1 \Delta)}^{-D} \nonumber\\
    & \leq \frac\epsilon2 f(\ve x^*, \ve z^*), \label{eq:estimate-gridopt-to-optimum}
  \end{align}
  where the last estimate is given by~\eqref{eq:lower-bound-for-optimum} in
  the case that $f$ is not constant on the feasible region.  On the other
  hand, if $f$ is constant, the
  estimate~\eqref{eq:estimate-gridopt-to-optimum} holds trivially.

  By \autoref{cor:gridproblem-approximation}
  we can compute a point $(\ve 
  x^m_{\epsilon/2}, \ve z^m_{\epsilon/2})\in P\cap\bigl(\frac1m\Z^{d_1}\times
  \Z^{d_2}\bigr)$ such
  that 
  \begin{equation}\label{eq:estimate-gridopt-to-gridapprox}
    (1-\tfrac\epsilon2) f(\ve x^m,\ve z^m) \leq  f(\ve x^m_{\epsilon/2}, \ve
    z^m_{\epsilon/2}) \leq f(\ve x^m, \ve z^m)
  \end{equation}
  in time polynomial in $\log m$, the encoding 
  length of $f$ and~$P$, the maximum total degree~$D$, and $1/\epsilon$. 
  Here $\log m$ is bounded by a polynomial in $\log M$, $D$ and $\log C$,
  so we can compute $(\ve x^m_{\epsilon/2}, \ve z^m_{\epsilon/2})$ 
  in time polynomial in the input
  size, the maximum total degree~$D$, and $1/\epsilon$.
  Now, using \eqref{eq:estimate-gridopt-to-gridapprox}
  and~\eqref{eq:estimate-gridopt-to-optimum}, we can estimate
  \begin{align*}
    f(\ve x^*, \ve z^*) - f(\ve x^m_{\epsilon/2}, \ve z^m_{\epsilon/2}) 
    & \leq f(\ve x^*,\ve z^*) - (1-\tfrac\epsilon2) f(\ve x^m,\ve z^m) \\
    & = \tfrac\epsilon2 f(\ve x^*,\ve z^*) + (1-\tfrac\epsilon2) \bigl( f(\ve
    x^*,\ve z^*) - f(\ve x^m, \ve z^m)\bigr) \\
    & \leq \tfrac\epsilon2 f(\ve x^*,\ve z^*) + \tfrac\epsilon2 f(\ve x^*,\ve z^*) \\
    & = \epsilon f(\ve x^*,\ve z^*). 
  \end{align*}
  Hence $f(\ve x^m_{\epsilon/2}, \ve z^m_{\epsilon/2}) \geq (1-\epsilon) f(\ve
  x^*, \ve z^*)$.\medbreak

  \noindent\emph{Part (b).}
  Let the dimension $d\geq2$ be fixed.  We prove that there does not exist a
  {\small PTAS} for the maximization of arbitrary polynomials over
  mixed-integer sets of polytopes. 
  We use the NP-complete problem AN1 on page 249 of \cite{GarJohn79}. This is
  to decide whether, given three positive integers $a,b,c$, there exists a
  positive integer $x<c$ such that $x^2 \equiv a\pmod b$.  This problem is
  equivalent to asking whether the maximum of the quartic polynomial
  function $f(x,y) = - (x^2-a-by)^2$ over the lattice points of the rectangle
  \begin{displaymath}
    P = \biggl\{\,(x,y) :  1 \leq x \leq c-1,\ \frac{1-a}{b} \leq y \leq
    \frac{(c-1)^2-a}{b} \,\biggr\}
  \end{displaymath}
  is zero or not.  If there existed a {\small PTAS} for the
  maximization of arbitrary polynomials over mixed-integer sets of polytopes,
  we could, for any fixed $0<\epsilon<1$, compute in polynomial time a solution
  $(x_\epsilon,y_\epsilon)\in P\cap\Z^2$ with $\bigl|f(x_\epsilon,y_\epsilon)
  - f(x^*,y^*)\bigr| \leq \epsilon \bigl|f(x^*,y^*)\bigr|$, where $(x^*,y^*)$ denotes an
  optimal solution.  Thus, we have $f(x_\epsilon,y_\epsilon)=0$ if and only if
  $f(x^*,y^*)=0$; this means we could solve the problem~AN1 in polynomial time.
\end{proof}

\section{Extension to arbitrary polynomials}
\label{section:weak}

In this section we drop the requirement of the polynomial being positive over the
feasible region.  As we showed
in \autoref{th:mipo-fptas}, 
there does not exist a {\small PTAS} for the
maximization of an arbitrary polynomial over polytopes in fixed dimension. 
We will instead show an approximation result like the one in
\cite{deklerk-laurent-parillo:ptas-polynomial-simplex}, i.e., we compute a
solution $(\ve x_\epsilon, \ve z_\epsilon)$ such that
\begin{equation}
  \label{eq:deklerk-approx}
  \bigl| f(\ve x_\epsilon, \ve z_\epsilon) - f(\ve x_{\max}, \ve z_{\max}) \bigr| 
  \leq \epsilon \bigl| f(\ve x_{\max}, \ve z_{\max}) - f(\ve x_{\min},\ve z_{\min}) \bigr|,
\end{equation}
where $(\ve x_{\max}, \ve z_{\max})$ is an optimal solution to the
maximization problem over the feasible region and $(\ve x_{\min}, \ve
z_{\min})$ is an optimal solution to the minimization problem.
Our algorithm has a running time that is polynomial in the input size, the
maximum total degree of~$f$, and $\frac1\epsilon$.  This means that while the result of
\cite{deklerk-laurent-parillo:ptas-polynomial-simplex} was a weak version of
a {\small PTAS} (for fixed degree), our result is a weak version of an
{\small FPTAS} (for fixed dimension).


The approximation algorithms for the integer case (\autoref{lemma:bounds}) and the
mixed-integer case (\autoref{th:mipo-fptas})
only work for polynomial objective functions that are non-negative on the
feasible region.  In order to apply them to an arbitrary polynomial objective
function~$f$, we need to add a constant term to $f$ that is large enough.  As
proposed in \cite{deloera-hemmecke-koeppe-weismantel:intpoly-fixeddim}, we can
use linear programming techniques to obtain a bound~$M$ on the variables and
then estimate 
\begin{displaymath}
  f(\ve x) \geq -rCM^D =: L_0,
\end{displaymath}
where $C$ is the largest absolute value
of a coefficient, $r$~is the number of monomials of~$f$, and $D$~is the
maximum total degree.  However, the range $\bigl| f(\ve x_{\max}, \ve
z_{\max}) - f(\ve x_{\min},\ve z_{\min}) \bigr|$ can be exponentially small
compared to $L_0$, so in order to obtain an approximation $(\ve x_\epsilon, \ve
z_\epsilon)$ satisfying~\eqref{eq:deklerk-approx}, we would need an
$(1-\epsilon')$-approximation to the problem of maximizing $g(\ve x,\ve z) :=
f(\ve x,\ve z) - L_0$ with an exponentially small value of~$\epsilon'$.

To address this difficulty, we will first apply an algorithm which will
compute an approximation $[L_i, U_i]$ of the range $[f(\ve x_{\min},\ve
z_{\min}),\allowbreak f(\ve x_{\max}, \ve z_{\max})]$ with constant quality.
To this end, we first prove a simple corollary of \autoref{th:mipo-fptas}.
\begin{corollary}[Computation of upper bounds for mixed-integer problems]
  \label{cor:mixed-dualbound}
  Let the dimension $d=d_1+d_2$ be fixed. 
  Let $P\subseteq\R^d$ be a
  rational convex polytope. 
  Let $f\in\Z[x_1,\dots,x_{d_1},\allowbreak z_1,\dots,z_{d_2}] 
  $ be a polynomial function with integer coefficients%
  and maximum total degree~$D$ that is
  non-negative on $P\cap \bigl(\R^{d_1}\times\Z^{d_2}\bigr)$.  Let
  $\delta>0$.  There exists an algorithm with running time polynomial in
  the input size, $D$, and~$\frac1{\delta}$ for computing an upper bound $u$ 
  such that
  \begin{equation}
    \label{eq:mixed-dualbound}
    f(\ve x_{\max},\ve z_{\max}) \leq u \leq (1+\delta) f(\ve x_{\max}, \ve z_{\max}),
  \end{equation}
  where $(\ve x_{\max},\ve z_{\max})$ is an optimal solution to the
  maximization problem of $f$ over $P\cap \bigl(\R^{d_1}\times\Z^{d_2}\bigr)$.
\end{corollary}
\begin{proof}
  Let $\epsilon = \frac{\delta}{1+\delta}$.  By
  \autoref{th:mipo-fptas}, we can,
  in time polynomial in the input size, $D$, and $\frac1\epsilon = 1 +
  \frac1{\delta}$, compute a
  solution $(\ve x_\epsilon, \ve z_\epsilon)$ with 
  \begin{equation}
    \bigl| f(\ve x_{\max}, \ve z_{\max}) - f(\ve x_\epsilon, \ve z_\epsilon) \bigr| 
    \leq \epsilon f(\ve x_{\max}, \ve z_{\max}).
  \end{equation}
  Let $u := \frac 1{1-\epsilon} f(\ve x_\epsilon, \ve z_\epsilon)  =
  (1+\delta) f(\ve x_\epsilon, \ve z_\epsilon)$.  
  Then
  \begin{equation}
    f(\ve x_{\max},\ve z_{\max}) \leq \frac1{1-\epsilon} f(\ve x_{\epsilon}, \ve
    z_{\epsilon}) = u 
  \end{equation}
  and
  \begin{align}
    (1+\delta) f(\ve x_{\max},\ve z_{\max}) 
    & \geq (1+\delta) f(\ve x_\epsilon, \ve z_\epsilon) \nonumber\\
    & = (1+\delta) (1-\epsilon) u \nonumber\\
    & = (1+\delta) \biggl(1-\frac{\delta}{1+\delta}\biggr) u
    = u.
  \end{align}
  This proves the estimate~\eqref{eq:mixed-dualbound}.
\end{proof}

\begin{algorithm}[Range approximation]\mbox{}\smallskip\par
  \label{algo:range}
  \noindent\textit{Input:} Mixed-integer polynomial optimization problem \eqref{eq:mixed-problem}, 
  a number $0<\delta<1$.
  
  \noindent\textit{Output:} Sequences $\{L_i\}$, $\{U_i\}$ of lower and upper
  bounds of~$f$ over the feasible region $P\cap
  \bigl(\R^{d_1}\times\Z^{d_2}\bigr)$ such that 
  \begin{equation}
    \label{eq:range}
    L_i \leq f(\ve x_{\min}, \ve z_{\min}) \leq f(\ve x_{\max}, \ve z_{\max}) \leq U_i
  \end{equation}
  and
  \begin{equation}
    \label{eq:range-limit}
    \lim_{i\to\infty} |U_i - L_i| = c \bigl(f(\ve x_{\max}, \ve z_{\max}) - f(\ve x_{\min}, \ve z_{\min})\bigr), 
  \end{equation}
  where $c$ depends only on the choice of~$\delta$.

  \begin{enumerate}
  \item By solving $2d$ linear programs over $P$, we find lower and upper
    integer bounds for each of the variables $x_1,\ldots,x_{d_1},z_1,\ldots,z_{d_2}$. Let $M$ be the
    maximum of the absolute values of these $2d$ numbers. Thus $|x_i|, |z_i|\leq M$
    for all $i$. Let $C$ be the maximum of the absolute values of all
    coefficients, and $r$ be the number of monomials of $f(x)$.  Then
    \[
    L_0:=-rCM^D \leq f(\ve x,\ve z) \leq rCM^D =: U_0,
    \]
    as we can bound the absolute value of each monomial of $f(x)$ by
    $CM^D$.
  \item Let $i := 0$.
  \item \label{algo:range:loop} 
    Using the algorithm of \autoref{cor:mixed-dualbound}, compute an upper
    bound~$u$  for the problem
    \begin{align*}
      \max\quad&g(\ve x,\ve z) := f(\ve x,\ve z) - L_i \\
      \st\quad &(\ve x,\ve z) \in P\cap \bigl(\R^{d_1}\times\Z^{d_2}\bigr)
    \end{align*}
    that gives
    a $(1+\delta)$-approximation to the optimal value.  Let $U_{i+1} := L_i + u$.
  \item Likewise, compute an upper
    bound~$u$ for the problem 
    \begin{align*}
      \max\quad & h(\ve x,\ve z) := U_i - f(\ve x,\ve z) \\
      \st\quad  & (\ve x,\ve z) \in P\cap \bigl(\R^{d_1}\times\Z^{d_2}\bigr)
    \end{align*}
    that gives a $(1+\delta)$-approximation to the optimal value.  
    Let $L_{i+1} := U_i - u$.
  \item $i := i+1$.
  \item Go to \ref{algo:range:loop}.
  \end{enumerate}
\end{algorithm}

\begin{lemma}\label{lemma:algo:range}
  \autoref{algo:range} is correct.  For fixed $0<\delta<1$, it computes
  the bounds $L_n$, $U_n$ satisfying
  \eqref{eq:range}~and~\eqref{eq:range-limit} in time polynomial in the input
  size and $n$. 
\end{lemma}
\begin{proof}
  We have
  \begin{equation}
    U_i - L_{i+1} \leq (1 + \delta) \bigl(U_i - f(\ve x_{\min}, \ve z_{\min})\bigr)
  \end{equation}
  and
  \begin{equation}
    U_{i+1} - L_i \leq (1 + \delta) \bigl(f(\ve x_{\max}, \ve z_{\max}) - L_i\bigr).
  \end{equation}
  This implies
  \begin{align*}
    U_{i+1} - L_{i+1} 
    &\leq \delta (U_i - L_i) + (1+\delta) \bigl(f(\ve x_{\max}, \ve z_{\max}) - f(\ve x_{\min}, \ve z_{\min})\bigr).
  \end{align*}
  Therefore
  \begin{align*}
    U_n - L_n &\leq \delta^n (U_0 - L_0) + (1 + \delta)
    \biggl(\sum_{i=0}^{n-2} \delta^i\biggr) \bigl(f(\ve x_{\max}, \ve
    z_{\max}) - f(\ve x_{\min}, \ve z_{\min})\bigr) \\
    & = \delta^n (U_0 - L_0) + (1 + \delta) \frac{1-\delta^{n-1}}{1-\delta}
    \bigl(f(\ve x_{\max}, \ve z_{\max}) - f(\ve x_{\min}, \ve
    z_{\min})\bigr)\\  
    & \to \frac{1+\delta}{1-\delta} \bigl(f(\ve x_{\max}, \ve z_{\max}) -
    f(\ve x_{\min}, \ve z_{\min})\bigr)\quad (n\to\infty). 
  \end{align*}

  The bound on the running time requires a careful analysis.  Because in each
  step the result $u$ (a rational number) of the bounding procedure
  (\autoref{cor:mixed-dualbound}) becomes part of the input in
  the next iteration, the encoding length of the input could grow exponentially
  after only polynomially many steps.  However, we will show that the encoding
  length only grows very slowly.

  First we need to remark that the auxiliary objective functions $g$~and~$h$ have
  integer coefficients except for the constant term, which may be rational.
  It turns out that the estimates in the proof of
  \autoref{th:mipo-fptas} (in
  particular, the local Lipschitz constant $L$ and the lower bound on the
  optimal value) are independent from the constant term of the objective
  function.  Therefore, the \emph{same} approximating grid $\frac 1m\Z^{d_1}\times
  \Z^{d_2}$ can be chosen in all iterations of \autoref{algo:range}; 
  the number $m$ only depends on $\delta$, the polytope~$P$, the maximum total
  degree~$D$, and the coefficients of~$f$ with the exception of the constant term.

  The construction in the proof of \autoref{cor:mixed-dualbound}
  obtains the upper bound~$u$ by multiplying the approximation $f(\ve x_\epsilon,\ve
  z_\epsilon)$ by $(1+\delta)$.  Therefore we have
  \begin{align}
    U_{i+1} & = L_i + u \nonumber\\
    & = L_i + (1+\delta) \bigl(f(\ve x_{\epsilon}, \ve z_{\epsilon}) - L_i\bigr)\nonumber\\
    & = -\delta L_i + (1+\delta) f(\ve x_{\epsilon}, \ve z_{\epsilon}).
    \label{eq:estimate-upper}
  \end{align}
  Because the solution $(\ve x_\epsilon,\ve
  z_\epsilon)$ lies in the grid $\frac 1m\Z^{d_1}\times \Z^{d_2}$, the value $f(\ve
  x_\epsilon, \ve z_\epsilon)$ is an integer multiple of $m^{-D}$.  
  This implies that, because $L_0 \leq f(\ve x_\epsilon, \ve z_\epsilon) \leq U_0$, the
  encoding length of the rational number $f(\ve x_\epsilon, \ve z_\epsilon)$
  is bounded by a polynomial in the input size of $f$~and~$P$. 
  Therefore the encoding length $U_{i+1}$ (and likewise $L_{i+1}$) only increases by an
  additive term that is bounded by a polynomial in the input size of $f$~and~$P$. 
\end{proof}


We are now in the position to prove \autoref{th:mipo-wfptas}.

\begin{proof}[Proof of \autoref{th:mipo-wfptas}]
  Clearly we can restrict ourselves to polynomials with integer coefficients.
  Let $m = (D+1) d_1 \Delta$, where $\Delta$ is the number from
  \autoref{lemma:grid-approximation-lowerdim}.  
  We apply \autoref{algo:range} using $0<\delta<1$ arbitrary to compute bounds $U_n$
  and $L_n$ for 
  \begin{displaymath}
    n = \bigl\lceil -\log_\delta \bigl(2 m^D (U_0 - L_0)\bigr)
    \bigr\rceil.
  \end{displaymath}
  Because $n$ is bounded by a polynomial in the input size and
  the maximum total degree $D$, this can be done in polynomial time.
  Now, by the proof of \autoref{lemma:algo:range}, we have 
  \begin{align}
    U_n - L_n &\leq \delta^n (U_0 - L_0) + (1 + \delta) \frac{1-\delta^{n-1}}{1-\delta}
    \bigl(f(\ve x_{\max}, \ve z_{\max}) - f(\ve x_{\min}, \ve z_{\min})\bigr)
    \nonumber \\
    & \leq \frac12 m^{-D} + \frac{1 + \delta}{1-\delta} \bigl(f(\ve x_{\max}, \ve
    z_{\max}) - f(\ve x_{\min},\ve z_{\min})\bigr).
  \end{align}
  If $f$ is constant on $P\cap\bigl( \R^{d_1}\times \Z^{d_2} \bigr)$, it
  is constant on $P\cap\bigl( \frac1m \Z^{d_1}\times \Z^{d_2} \bigr)$, 
  then $U_n - L_n \leq\frac12 m^{-D}$.  Otherwise,
  by \autoref{lemma:mixed-constantness}, we have
  $U_n - L_n \geq f(\ve 
  x_{\max}, \ve z_{\max}) - f(\ve x_{\min}, \ve z_{\min}) \geq m^{-D}$.  This
  settles part~(a). 

  For part (b), if $f$ is constant on $P\cap\bigl( \R^{d_1}\times \Z^{d_2}
  \bigr)$, we return an arbitrary solution as an 
  optimal solution.  Otherwise, we can estimate further:
  \begin{align}\label{eq:mixed-range-estimate}
    U_n - L_n &\leq \biggl(\frac12 + \frac{1 + \delta}{1-\delta}\biggr)
    \bigl(f(\ve x_{\max}, \ve z_{\max}) - f(\ve x_{\min}, \ve z_{\min})\bigr).
  \end{align}
  
  Now we apply the algorithm of \autoref{th:mipo-fptas} to the 
  maximization problem of the polynomial function $f' := f - L_n$, which is
  non-negative over the feasible region $P\cap \bigl(\R^{d_1}\times \Z^{d_2}\bigr)$. 
  We compute a point $(\ve x_{\epsilon'}, \ve z_{\epsilon'})$ where
  $\epsilon' = \epsilon {\bigl(\frac12 + \frac{1 + \delta}{1-\delta}\bigr)}^{-1} $ such that 
  \begin{displaymath}
    \bigl| f'(\ve x_{\epsilon'}, \ve z_{\epsilon'}) - f'(\ve x_{\max}, \ve z_{\max}) \bigr| \leq \epsilon' f'(\ve x_{\max}, \ve z_{\max}).
  \end{displaymath}
  Then we obtain the estimate
  \begin{align*}
    \bigl| f(\ve x_{\epsilon'}, \ve z_{\epsilon'}) - f(\ve x_{\max}, \ve z_{\max}) \bigr| 
    &\leq \epsilon' \bigl( f(\ve x_{\max}, \ve z_{\max}) - L_n \bigr) \\
    &\leq \epsilon' \bigl( U_n - L_n \bigr) \\
    &\leq \epsilon' \biggl(\frac12 + \frac{1 + \delta}{1-\delta}\biggr)
    \bigl(f(\ve x_{\max}, \ve z_{\max}) - f(\ve x_{\min}, \ve z_{\min})\bigr) \\
    & = \epsilon \bigl(f(\ve x_{\max},\ve z_{\max}) - f(\ve x_{\min}, \ve z_{\max})\bigr),
  \end{align*}
  which proves part (b).
\end{proof}


\bibliographystyle{amsplain}
\bibliography{barvinok,iba-bib,weismantel}

\providecommand\FFull[1]{#1} \providecommand\Full[1]. \providecommand\NP{NP}
  \providecommand\textem[1]{{\textbf{#1}}}
  \providecommand\with{{\kern1pt--\kern1pt}}
  \def\polhk#1{\setbox0=\hbox{#1}{\ooalign{\hidewidth
  \lower1.5ex\hbox{`}\hidewidth\crcr\unhbox0}}}
\providecommand{\bysame}{\leavevmode\hbox to3em{\hrulefill}\thinspace}
\providecommand{\MR}{\relax\ifhmode\unskip\space\fi MR }
\providecommand{\MRhref}[2]{%
  \href{http://www.ams.org/mathscinet-getitem?mr=#1}{#2}
}
\providecommand{\href}[2]{#2}
\begin{thebibliography}{10}

\bibitem{bar}
Alexander~I. Barvinok, \emph{Polynomial time algorithm for counting integral
  points in polyhedra when the dimension is fixed}, Mathematics of Operations
  Research \textbf{19} (1994), 769--779.

\bibitem{BarviPom}
Alexander~I. Barvinok and James~E. Pommersheim, \emph{An algorithmic theory of
  lattice points in polyhedra}, New Perspectives in Algebraic Combinatorics,
  Math. Sci. Res. Inst. Publ., vol.~38, Cambridge Univ. Press, Cambridge, 1999,
  pp.~91--147.

\bibitem{bellare-rogaway-1993}
Mihir Bellare and Philip Rogaway, \emph{The complexity of aproximating a
  nonlinear program}, in Pardalos
  \cite{pardalos:complexity-numerical-optimization}.

\bibitem{CoxLittleOShea92}
David~A. Cox, John~B. Little, and Donal O'Shea, \emph{Ideals, varieties, and
  algorithms: An introduction to computational algebraic geometry and
  commutative algebra}, Springer, Berlin, Germany, 1992.

\bibitem{deklerk-laurent-parillo:ptas-polynomial-simplex}
Etienne de~Klerk, Monique Laurent, and Pablo~A. Parrilo, \emph{A {PTAS} for the
  minimization of polynomials of fixed degree over the simplex}, Theoretical
  Computer Science \textbf{361} (2006), 210--225.

\bibitem{latte2}
Jes{\'u}s~A. De~Loera, David Haws, Raymond Hemmecke, Peter Huggins, Bernd
  Sturmfels, and Ruriko Yoshida, \emph{Short rational functions for toric
  algebra and applications}, Journal of Symbolic Computation \textbf{38}
  (2004), no.~2, 959--973.

\bibitem{deloera-hemmecke-koeppe-weismantel:intpoly-fixeddim}
Jes{\'u}s~A. De~Loera, Raymond Hemmecke, Matthias K{\"o}ppe, and Robert
  Weismantel, \emph{Integer polynomial optimization in fixed dimension},
  Mathematics of Operations Research \textbf{31} (2006), no.~1, 147--153.

\bibitem{latte1}
Jes{\'u}s~A. De~Loera, Raymond Hemmecke, Jeremiah Tauzer, and Ruriko Yoshida,
  \emph{Effective lattice point counting in rational convex polytopes}, Journal
  of Symbolic Computation \textbf{38} (2004), no.~4, 1273--1302.

\bibitem{GarJohn79}
Michael~R. Garey and David~S. Johnson, \emph{Computers and intractability: A
  guide to the theory of {NP}-completeness}, Freeman, San Francisco, 1979.

\bibitem{Hastad:inapprox97}
Johan H{\aa}stad, \emph{Some optimal inapproximability results}, Proceedings of
  the 29th Symposium on the Theory of Computing (STOC), ACM, 1997, pp.~1--10.

\bibitem{jones-1982}
James~P. Jones, \emph{Universal diophantine equation}, Journal of Symbolic
  Logic \textbf{47} (1982), no.~3, 403--410.

\bibitem{Len83}
Hendrik~W. Lenstra, Jr., \emph{Integer programming with a fixed number of
  variables}, Mathematics of Operations Research \textbf{8} (1983), 538--548.

\bibitem{matiyasevich-1970}
{Yu}ri~V. Matiyasevich, \emph{Enumerable sets are diophantine}, Doklady
  Akademii Nauk SSSR \textbf{191} (1970), 279--282, (Russian); English
  translation, Soviet Mathematics Doklady, vol. 11 (1970), pp. 354--357.

\bibitem{matiyasevich-1993}
\bysame, \emph{Hilbert's tenth problem}, The MIT Press, Cambridge, MA, USA,
  1993.

\bibitem{pardalos:complexity-numerical-optimization}
Panos~M. Pardalos (ed.), \emph{Complexity in numerical optimization}, World
  Scientific, 1993.

\bibitem{Renegar:1992:CCGa}
James Renegar, \emph{On the computational complexity and geometry of the
  first-order theory of the reals, part {I}: Introduction. {Preliminaries}.
  {The} geometry of semi-algebraic sets. {The} decision problem for the
  existential theory of the reals}, Journal of Symbolic Computation \textbf{13}
  (1992), no.~3, 255--300.

\bibitem{Renegar:1992:CCGb}
\bysame, \emph{On the computational complexity and geometry of the first-order
  theory of the reals, part {II}: The general decision problem. {Preliminaries}
  for quantifier elimination}, Journal of Symbolic Computation \textbf{13}
  (1992), no.~3, 301--328.

\bibitem{Renegar:1992:CCGc}
\bysame, \emph{On the computational complexity and geometry of the first-order
  theory of the reals. part {III}: Quantifier elimination}, Journal of Symbolic
  Computation \textbf{13} (1992), no.~3, 329--352.

\bibitem{Renegar:1992:Approximating}
\bysame, \emph{On the computational complexity of approximating solutions for
  real algebraic formulae}, SIAM Journal on Computing \textbf{21} (1992),
  no.~6, 1008--1025.

\bibitem{vavasis-1993}
Stephen~A. Vavasis, \emph{Polynomial time weak approximation algorithms for
  quadratic programming}, in Pardalos
  \cite{pardalos:complexity-numerical-optimization}.

\end{thebibliography}

\end{document}